\documentclass[12pt]{article}

\usepackage[latin1]{inputenc}
\usepackage{amsmath,amssymb}
\usepackage{latexsym}

\newtheorem{theorem}{Theorem}[section]

\newtheorem{proposition}[theorem]{Proposition}

\newtheorem{remark}[theorem]{Remark}
\numberwithin{equation}{section}

\def\proof{{\medskip\noindent {\bf Proof. }}}

\def\qed{{\hfill $\square$ \bigskip}}

\def\square{{\vcenter{\vbox{\hrule height.3pt
        \hbox{\vrule width.3pt height5pt \kern5pt
           \vrule width.3pt}
        \hrule height.3pt}}}}

  \def\sC {{\cal C}}

  \def\sL {{\cal L}}
\def\sM {{\cal M}}  
  
\def\sS {{\cal S}}

\def\E{{\mathbb E}}
\def\P{{\mathbb P}}
\def\norm#1{{\Vert #1 \Vert}}
\def\del{{\partial}}
\def\lam{{\lambda}}

\def\bea{\begin{align*}}
\def\eea{\end{align*}}
\def\bee{\begin{equation}}
\def\eee{\end{equation}}

\def\R{{\mathbb R}}
\def\E{{{\mathbb E}\,}}
\def\P{{\mathbb P}}

\def\lam{{\lambda}}

\def\al{{\alpha}}
\def\grad{{\nabla}}
\def\proof{{\medskip\noindent {\bf Proof. }}}

\def\qed{{\hfill $\square$ \bigskip}}

\def\eps{\varepsilon}

\def\norm#1{\Vert #1 \Vert}

 \def\qq {\qquad}
\def\del{{\partial}}

\def\bs{\bigskip}
\def\cl#1{\centerline{#1}}

\def\square{{\vcenter{\vbox{\hrule height.3pt
        \hbox{\vrule width.3pt height5pt \kern5pt
           \vrule width.3pt}
        \hrule height.3pt}}}}

\def\tfrac#1#2{{\textstyle {\frac{#1}{#2}}}}

\def\tlint{{- \kern-0.85em \int \kern-0.2em}}  
\def\dlint{{- \kern-1.05em \int \kern-0.4em}}  

  \def\sC {{\cal C}}

  \def\sL {{\cal L}}
\def\sM {{\cal M}}  
  
\def\sS {{\cal S}}

\def\nn{{\nonumber}}

\begin{document}

\title{A new technique for proving uniqueness for martingale problems}

\author{Richard F. Bass\footnote{Research partially supported by NSF grant
DMS-0244737.}  \; and Edwin Perkins\footnote{Research partially supported
by NSERC.}}

\date{\today}

\maketitle

\cl{Version:{Bismut paper; c.tex}}\bs

\begin{abstract}  
\noindent A new technique for proving uniqueness of martingale
problems is introduced. The method is illustrated in the context
of elliptic diffusions in $\R^d$.

%
\end{abstract}

\section{Introduction}

When trying to prove uniqueness of a stochastic process corresponding
to an operator, one of the
most useful approaches is to consider the associated martingale problem.
If $\sL$ is an operator and $w$ is a point in the state space $\sS$, a probability $\P$ on the set of paths $t\to X_t$ taking values in $\sS$ is a solution of the
martingale problem for $\sL$ started at $w$ if $\P(X_0=w)=1$ and $f(X_t)-f(X_0)-\int_0^t \sL f(X_s)\, ds$
is a martingale with respect to $\P$ for every $f$ in an appropriate class $\sC$
of functions.

The archetypical example is to let 
\bee\label{defL}
\sL f(x)=\sum_{i,j=1}^d a_{ij}(x) D_{ij}f(x).
\eee
Here, and for the rest of this paper, the state space $\sS$ is $\R^d$, the probability measure is on the set
of functions that are continuous maps from $[0,\infty)$ into
$\R^d$ with the $\sigma$-field generated by the cylindrical sets,
$D_{ij}f =\del^2 f/\del x_i \del x_j$, and the class $\sC$ of functions
is the collection $C^2_b$ of $C^2$ functions which are bounded and whose first and
second partial derivatives are bounded. 

Stroock and Varadhan introduced the notion of martingale problem and
proved in the case above that there was existence and uniqueness of
the solution to the martingale problem provided the $a_{ij}$ were
bounded and continuous in $x$ and the matrix $a(x)$ was strictly positive definite
for each $x$. See \cite{B2} or \cite{SV} for an account of this result.

In this paper we present a new method of proving uniqueness for 
martingale problems. We illustrate it for the operator $\sL$ given
in (\ref{defL}) under the assumption that the $a_{ij}$ are H\"older continuous
in $x$. Our proof does not give as strong a result as that of Stroock and
Varadhan in that we require H\"older continuity. (Actually, we only require
a Dini-like condition, but this is still more than just
requiring continuity.) In fact, when the $a_{ij}$ are H\"older continuous,
an older method using Schauder estimates can be applied.

Nevertheless our technique is applicable to situations for which no other
known method seems to work. A precursor of our method, much disguised,  was used in 
\cite{B1} to prove uniqueness for pure jump processes which were of
variable order, i.e., the operator can not be viewed as a perturbation
of a symmetric stable process of any fixed order. The result of
\cite{B1} was improved in \cite{T} to allow more general jump processes.
Moreover our technique is useful in problems arising from
certain infinite dimensional situations in the theory of
stochastic partial differential equations and the theory of superprocesses; 
see \cite{BP}. Finally, even in the elliptic diffusion case considered here, the proof is 
elementary and short.

Stroock and Varadhan's method was essentially to view $\sL$ given
in (\ref{defL}) as a perturbation of the Laplacian with respect to
the space $L^p$ for appropriate $p$. 
The method using Schauder estimates views $\sL$
as a perturbation of the Laplacian with respect to the H\"older space
$C^\al$ for appropriate $\al$. We use a quite different approach. We
view $\sL$ as a mixture of constant coefficient operators and 
use a mixture of the corresponding semigroups as an approximation
of the semigroup for $\sL$.

We use our method to prove the following theorem.

\begin{theorem}\label{main}
Suppose $\sL$ is given by (\ref{defL}), the matrices $a(x)$ are
bounded and uniformly positive definite, and there exist $c_1$ and
$\al$ such that 
\bee\label{Holder}
|a_{ij}(x)-a_{ij}(y)|\leq c_1(1\land |x-y|^\al)
\eee
for all $i,j=1, \ldots, d$ and all $x,y\in \R^d$. Then for each $w\in \R^d$
the
solution to the martingale problem for $\sL$ started at $w$ is unique.
\end{theorem}

We do not consider existence, since that is much easier, and we have
nothing to add to the existing proofs. The same comment
applies to the inclusion of drift terms. In Section 2 we  give some
easy estimates and in Section 3 we prove Theorem \ref{main}.
The letter $c$ denotes constants whose exact value is unimportant and
may change from occurrence to occurrence.

\section{Some estimates}

All the matrices we consider will be $d$ by $d$, bounded, symmetric, and uniformly elliptic,
that is, there exist constants $\Lambda_m$ and $\Lambda_M$ such
that
\bee\label{elliptic}
\Lambda_m\sum_{i=1}^d z_i^2\leq \sum_{i,j=1}^d a_{ij}z_iz_j\leq \Lambda_M\sum_{i=1}
^d z_i^2, \qq (z_1, \ldots, z_d)\in \R^d.
\eee
Given any such matrix $a$, we use $A$ for $a^{-1}$. It follows easily that 
\bee\label{bnd}
\sup_j\Big(\sum_{i=1}^da_{ij}^2\Big)^{1/2}\le \Lambda _M,\quad \sup_j\Big(\sum_{i=1}^d A_{ij}^2\Big)^{1/2}\le \Lambda_m^{-1}
\eee

Define 
\bee\label{ptdef}
p^a(t,x,y)=(2\pi t)^{-d/2} (\det a)^{-1/2} e^{-(y-x)^TA(y-x)/(2t)},
\eee
and let
\bee\label{Ptdef}
P_t^a f(x)=\int p^a(t,x,y)f(y)\, dy
\eee
be the corresponding transition operator.
We assume throughout that the matrix valued function $a(y)$ satisfies 
the hypotheses of Theorem~\ref{main} and \eqref{elliptic}.  
Note that for $a$ fixed, $p^a(t,x,y)\, dy$ is a Gaussian distribution for
each $x$, but that $p^{a(y)}(t,x,y)\, dy$ need not be a probability
measure.
All numbered constants will depend only $\Lambda_m,\Lambda_M$ and $d$.  

We have the following. 

\begin{proposition}\label{P2.1}
There exist $c_1$, $c_2$ and a function $c_3(p),p>0$, depending only on $\Lambda_M$ and $\Lambda_m$, such that for all $t,N,p>0$ and $x\in\R^d$, 

(a) $\int p^{a(y)}(t,x,y)\, dy\leq c_1$.

(b) $$\int_{|y-x|>N/\sqrt t} p^{a(y)}(t,x,y)\, dy\leq c_1e^{-c_2N^2}.$$

(c) 
For each $i\le d$, 
$$\int \Big(\frac{|x_i-y_i|^2}{t}\Big)^p p^{a(y)}(t,x,y)\, dy\leq c_3(p).$$

\end{proposition}

\proof For (a), after a change of variables $z=(y-x)/\sqrt t$, we need to bound
\begin{align*}
\int (2\pi)^{-d/2}&(\det a(x+z\sqrt t))^{-1/2}
e^{-z^T A(x+z\sqrt t)z/2} \, dz\\
&\leq \Big(\frac{\Lambda_M}{\Lambda_m}\Big)^{d/2}
\int (2\pi \Lambda_M)^{-d/2} e^{-z^Tz/2\Lambda_M}\, dz
\leq  \Big(\frac{\Lambda_M}{\Lambda_m}\Big)^{d/2}.
\end{align*}
(b) and (c) are similar.
\qed

Let $\norm{f}$ be the $C_0$ norm of $f$.

\begin{proposition}\label{P2.2}
Let $g\in C^2$ with compact support and let
$$F_\eps(x)=\int g(y)p^{a(y)}(\eps^2,x,y)\, dy.$$
Then
$F_\eps(x)$ converges to $g(x)$ boundedly and pointwise.
\end{proposition}

\proof 
Because $g$ is bounded, using  Proposition \ref{P2.1}(a) we see that
the quantity
$\sup_{\eps>0}\Vert F_\eps\Vert$ is finite. 
We next consider pointwise convergence.
After a change of variables, we have
$$F_\eps(x)=\int g(x+\eps z) (2\pi)^{-d/2} (\det a(x+\eps z))^{-1/2} e^{-z^TA(x+\eps z)z/2} \, dz.$$
Since $|g(x+\eps z)-g(x)|\leq \eps |z|\, \norm{\grad g}$, $F_\eps$
differs from
$$g(x) \int (2\pi)^{-d/2} (\det(a(x+\eps z))^{-1/2} e^{-z^TA(x+\eps z)z/2} \, dz$$
by at most 
$$\norm{\grad g}\int (2\pi)^{-d/2}(\det(a(x+\eps z)))^{-1/2} \eps |z|e^{-z^TA(x+\eps z)z/2} \, dz,$$
and this goes to 0 as $\eps\to 0$ by a change of variables and Proposition~\ref{P2.1}(c) with $p=1/2$. 
Let $$V(\eps,x,z)=(2\pi)^{-d/2}(\det(a(x+\eps z)))^{-1/2} e^{-z^TA(x+\eps z)z/2}.$$
It therefore suffices to show $$\int V(\eps,x,z)\, dz\to
\int V(0,x,z)\, dz,$$
where we note this right-hand side is $1$.
  Using Proposition \ref{P2.1}(b) and the same change of variables, it suffices to show
$$\int_{|z|\leq N} V(\eps,x,z)\,dz\to \int_{|z|\leq N} V(0,x,z)\, dz.$$
But this last follows by dominated convergence.
\qed

\begin{proposition}\label{P2.3}
There exists a constant $c_4$ such that
\begin{equation*}
\int |a_{ij}(y)-a_{ij}(x)|\, |D_{ij}p^{a(y)}(t,x,y)|\, dy
\leq \begin{cases} c_4t^{\frac{\al}{2}-1}, & t\leq 1,\\
 c_4t^{-1},& t\geq 1.
\end{cases}
\end{equation*}
\end{proposition}
\proof  A computation shows that
\begin{align}
D_{ij}&p^{a(y)}(t,x,y) \label{Ed.3.1}\\ 
&=t^{-1}p^{a(y)}(t,x,y)\Bigl[\sum_k\sum_l{(y_k-x_k)A_{ki}(y)A_{lj}(y)(y_l-x_l)\over t}-A_{ij}(y)\Bigr].\nn
\end{align}
By \eqref{bnd} and Cauchy-Schwarz we have 
\begin{align}
\int& |a_{ij}(y)-a_{ij}(x)|\, |D_{ij}p^{a(y)}(t,x,y)|\, dy\nonumber\\
&\leq \Bigl[\int|a_{ij}(y)-a_{ij}(x)|t^{-1}p^{a(y)}(t,x,y)[|x-y|^2t^{-1}\Lambda_m^{-2}+\Lambda_m^{-1}]\,dy.\label{genbnd}
\end{align}
Suppose first that $t\le 1$.  By the H\"older condition on $a$ the above is at most
\begin{align*}&c\int{|y-x|^\alpha\over t^{\alpha/2}}\Bigl[{|x-y|^2\over t}+1\Bigr]p^{a(y)}(t,x,y)\,dy\,t^{\alpha/2-1}\\
&\le ct^{\alpha/2-1},
\end{align*}
where we have used Proposition~\ref{P2.1}(c) in the last inequality. 

For the case $t>1$ simply use the boundedness of $a$ in \eqref{genbnd} and Proposition~\ref{P2.1} again to bound it by $ct^{-1}$.

\section{Proof of Theorem \ref{main}}

For $f\in C_b^2$ and $a$ a matrix with constant coefficients define
$$\sM^a f(x)=\sum_{i,j=1}^d a_{ij} D_{ij}f(x).$$
Define the corresponding semigroup by (\ref{Ptdef}),
and let
$R^a_\lam f=\int_0^\infty e^{-\lam t} P_t^a f\, dt.$
For $f\in C_b^2$ we have
$$\sL f(x)=\sM^{a(x)} f(x).$$
Note that
\bee\label{E3.1}(\lam-\sM^{a(y)})R_\lam^{a(y)} P_\eps^{a(y)} f(x)
=P_\eps^{a(y)} f(x).
\eee
One way to verify that the superscript $a(y)$ does not cause 
any difficulty here  is to check that
$$\sum_{i,j=1}^d a_{ij}(y) \frac{\del^2}{\del x_i \del x_j}
p^{a(y)}(s,x,y)=\frac{\del}{\del s} p^{a(y)}(s,x,y),$$
and then in the definition of $R_\lam^{a(y)}$ use integration by parts in the time variable.
By replacing $\eps$ with $\eps/2$, setting $f(z)=p^{a(y)}(\eps/2,z,y)$ and using Chapman-Kolmogorov, we see that \eqref{E3.1} implies
\bee\label{E3.1td}(\lam-\sM^{a(y)})(R_\lam^{a(y)} p^{a(y)}(\eps,\cdot,y))(x)=p^{a(y)}(\eps,x,y).
\eee
We are now ready to prove Theorem \ref{main}.

\proof Suppose $\P_1, \P_2$ are two solutions to the martingale
problem for $\sL$ started at a point $w$. Define
$$S_\lam^i f=\E_i\int_0^\infty e^{-\lam t} f(X_t)\, dt, \qq i=1,2,$$
and $$S_\lam^\Delta f=S_\lam^1 f- S_\lam^2.$$
We make two observations. First, because $\P_i$ need not come from a 
Markov process, $S^i_\lam f$ is not a function, and so $S_\lam^\Delta$ is a linear
functional. Second, if $$\Theta=\sup_{\norm{f}\leq 1} |S_\lam^\Delta f|,$$
then $\Theta<\infty$.

If $f\in C_b^2$, then by the definition of the martingale problem
$$\E_i f(X_t)-f(w)=\E_i\int_0^t \sL f(X_s)\, ds, \qq i=1,2.$$
Multiply both sides by $\lam e^{-\lam t}$, integrate over
$t$ from 0 to $\infty$, and use Fubini to obtain
$$f(w)=S_\lam^i(\lam f-\sL f), \qq i=1,2,$$
or
\bee S_\lam^\Delta (\lam f-\sL f)=0.\label{resz}\eee

Let $g\in C^2$ with compact support and set
$$f_\eps(x)=\int R_\lam^{a(y)} (p^{a(y)}(\eps,\cdot,y))(x)g(y)\, dy.$$
Since this is the same as
$$e^{-\lam \eps} \int \int_\eps^\infty e^{-\lam t} p^{a(y)}(t,x,y) \, dt\, g(y)\, dy,$$
we see that $f_\eps$ is in $C_b^2$ in $x$ by dominated convergence.

To calculate $(\lam-\sL)f_\eps$ it is easy to differentiate under the $dy$ integral and so we may write
\begin{align*}
(\lam -\sL)f_\eps(x)&=(\lam-\sM^{a(x)})f_\eps(x)\\
&=\int (\lam-\sM^{a(y)})R_\lam^{a(y)}(p^{a(y)}(\eps,\cdot,y))(x)g(y)\, dy\\
&~~~+\int (\sM^{a(y)}-\sM^{a(x)})R_\lam^{a(y)}(p^{a(y)}(\eps,\cdot,y))(x)g(y)\, dy\\
&:= I_\eps(x)+J_\eps(x).
\end{align*}
 By Proposition \ref{P2.3}, 
\begin{align*}|J_\eps(x)|&\leq \sum_{i,j=1}^d\int_0^\infty e^{-\lambda t}\int|a_{ij}(y)-a_{ij}(x)|\\
&\qq \qq\times |D_{ij}p^{a(y)}(\eps+t,x,y)|\,|g(y)|\,dy\,dt\\
&\le d^2\Vert g\Vert\int_0^\infty e^{-\lambda t}c_4t^{-1}(t^{\alpha/2}\wedge 1)\,dt\cr
&\le {1\over 2}\Vert g\Vert,
\end{align*}
for $\lam\ge \lam_0(\alpha,d,c_4)$.  
By (\ref{E3.1td}), $I_\eps(x)=\int p^{a(y)}(\eps,x,y)g(y)\, dy$, and so by
Proposition \ref{P2.2}, $I_\eps(x)$ converges to $g$ boundedly and
pointwise.
Since $S_\lam^\Delta (\lam-\sL) f_\eps=0$ by \eqref{resz}, we have
$|S_\lam^\Delta I_\eps|=|S_\lam^\Delta J_\eps|$. Letting $\eps\to 0$,
$$|S_\lam^\Delta g|=\lim_{\eps\to 0}  |S_\lam^\Delta I_\eps|=\lim_{\eps\to 0} |S_\lam^\Delta J_\eps|
\leq \Theta \limsup_{\eps\to 0} \norm{J_\eps}\leq \tfrac12 \Theta \norm{g}.$$
Using a monotone class argument, the above inequality holds for all
bounded $g$, and then taking the supremum over $g$ such that
$\norm{g}\leq 1$, we have $\Theta\leq \frac12 \Theta$. Since $\Theta<
\infty$, this implies that $\Theta=0$.

From this point on, we use standard arguments. 
By the uniqueness of the Laplace transform together
with continuity in $t$, $\E_1 f(X_t)=\E_2 f(X_t)$ for all $t$ if
$f$ is continuous and bounded. Using regular conditional probabilities,
one shows as usual that the finite dimensional distributions under
$\P_1$ and $\P_2$ agree. This suffices to prove uniqueness;
see \cite{B2} or \cite{SV} for details. 
\qed

Note that no localization argument is needed in the above proof.

\medskip

\begin{minipage}[t]{0.39\textwidth}
{\bf Richard F. Bass}\\
Department of Mathematics\\
University of
Connecticut \\
Storrs, CT 06269-3009, USA\\
{\it bass@math.uconn.edu}
\end{minipage}
\hfill
\begin{minipage}[t]{0.55\textwidth}
{\bf Edwin Perkins}\\
Department of Mathematics\\
University of British Columbia\\
Vancouver, B.C. V6T 1Z2, Canada\\
{\it perkins@math.ubc.ca} 
\end{minipage}

\end{document}